\documentclass[twoside,leqno,10pt]{amsart}
\usepackage{amsfonts}
\usepackage{amsmath}
\usepackage{amscd}
\usepackage{amssymb}
\usepackage{amsthm}
\usepackage{amsrefs}
\usepackage{latexsym}
\usepackage{bbm}
\numberwithin{equation}{section}
\begin{document}
\newtheorem{theorem}{Theorem}
\newtheorem{lemma}{Lemma}
\newtheorem{prop}{Proposition}
\newtheorem{corollary}{Corollary}
\newtheorem{conjecture}{Conjecture}
\numberwithin{equation}{section}
\newcommand{\dif}{\mathrm{d}}
\newcommand{\intz}{\mathbb{Z}}
\newcommand{\ratq}{\mathbb{Q}}
\newcommand{\natn}{\mathbb{N}}
\newcommand{\comc}{\mathbb{C}}
\newcommand{\rear}{\mathbb{R}} 
\newcommand{\prip}{\mathbb{P}}
\newcommand{\uph}{\mathbb{H}}

\title{Bombieri-Vinogradov Type Theorem for Sparse Sets of Moduli}
\author{Stephan Baier \and Liangyi Zhao}
\date{\today}
\maketitle

\begin{abstract}
In this paper, we establish theorems of Bombieri-Vinogradov type and Barban-Davenport-Halberstam type for sparse sets of moduli.  As an application, we prove that there exist infinitely many primes of the form $p=am^2+1$ such that $a\le p^{5/9+\varepsilon}$.
\end{abstract}



\section{Introduction and History}

The classical Bombieri-Vinogradov theorem \cites{EB1, AIV} touches upon the distribution of primes in arithmetic progressions on average.  More precisely, the theorem states the following. \newline

Let $A>0$ be fixed.  Then
\begin{equation} \label{classicbomvino}
\sum_{q \leq x^{1/2}/(\log x)^{A+5}} \max_{y\leq x} \max_{\substack{a \\ \gcd(a,q)=1}} \left| \psi(y;q,a) - \frac{y}{\varphi(q)} \right| \ll \frac{x}{(\log x)^A}.
\end{equation}
Here and after, $\varphi(n)$ is the Euler function and 
\[ \psi(x;q,a) = \sum_{\substack{n \leq x \\ n \equiv a \bmod{q} }} \Lambda(n), \]
where $\Lambda(n)$ is the von Mangoldt function. \newline

Moreover, theorems concerning the mean square error in the prime number theorem for arithmetic progressions were initiated by Barban \cite{MBB2} and Davenport and Halberstam \cites{HDHH1}.  Their results were sharpened by Gallagher \cite{PXG}.  The theorem in question states the following.  For any $A>0$ fixed, we have
\begin{equation} \label{classicbardavhal}
\sum_{q\leq x/(\log x)^{A+1}} \sum_{\substack{a=1 \\ \gcd(a,q)=1}}^{q} \left( \psi(x; q,a) - \frac{x}{\varphi(q)} \right)^2 \ll \frac{x^2}{(\log x)^A}.
\end{equation}
Montgomery \cite{HM4} showed that \eqref{classicbardavhal} is the best possible. \newline

In this paper, we aim to prove theorems of the types of Bombieri-Vinogradov and Barban-Halberstam-Davenport for sparse sets of moduli.  Theorems of both types are also proved for square moduli.  As in the classical theorems, one of the key ideas in the proof is the large sieve.  We shall need versions of the large sieve for these sparse sets of moduli.  Works in the direction of large sieve for special moduli can be found in \cites{Ba1, Zha, LZ4, SBLZ, SBLZ2}. \newline

It is also note-worthy that the implied constants in both \eqref{classicbomvino} and \eqref{classicbardavhal} are ineffective because of the use of the Siegel-Walfisz theorem in the estimates.  Such ineffectiveness is also present in Theorem~\ref{genBDHsquare} and Theorem~\ref{bomvinosquare} of this paper for the same reason.

\section{Statements of Results}

Throughout this paper, we assume without loss of generality that $x\ge 3$. 
We start with the following notation and definition.  Let $S$ be a set of natural numbers.  For $t\in \natn$ we set
\begin{equation} \label{defwelldis}
S_t = \{ q \in \natn : qt \in S \} , \; S_t(R) = \{ q \in S_t : R < q \leq 2R \}
\end{equation}
and
\[ S(R) = S_1(R). \]
We say that $S$ is well-distributed if, for $\varepsilon>0$, $t\in \natn$, $R \leq x < x+y \leq 2R$ and $\gcd ( k,l) =1$ we have
\begin{equation} \label{expected}
\left| \left\{ q \in S_t : x \leq q \leq x+y, \; q \equiv l \pmod{k} \right\} \right| \ll \left( \frac{|S_t(R)|y}{kR} + 1 \right) (Rt)^{\varepsilon},
\end{equation}
where the implied constant depends at most on $\varepsilon$ and $S$.  It should be noted that the above is a natural definition of a set being ``well-distributed'' in arithmetic progressions as the majorant in \eqref{expected} gives essentially the ``expected'' cardinality of the left-hand side of \eqref{expected}.  \newline

We are now ready to state the results of this paper.  First we have the following Barban-Davenport-Halberstam type theorem for sparse sets of moduli.

\begin{theorem} \label{genBDH}
Let $S$ be a well-distributed set.  Assume that
\begin{equation} \label{conditionforboth}
\left| S_q \left( \frac{Q}{q} \right) \right| \ll \frac{|S(Q)|}{q^{\varepsilon}}
\end{equation}
for any fixed $\varepsilon >0$ and all $q \in \natn$ with the implied constant depending at most on $\varepsilon$ and $S$.  Suppose further that
\begin{equation} \label{conditionfor1}
 |S(Q)| \gg \sqrt{Q}, \; \mbox{as} \; Q \to \infty.
\end{equation}
Then if $Q \leq x^{1-\varepsilon}$, we have, for any $A>0$ fixed,
\begin{equation} \label{genBDHeq}
\sum_{q\in S(Q)} \sum_{\substack{a=1 \\ \gcd(a,q)=1}}^q \left| \psi(x;q,a) - \frac{x}{\varphi(q)} \right|^2 \ll \frac{|S(Q)|}{Q} \frac{x^2}{(\log x)^A}.
\end{equation}
\end{theorem}

As a corollary to Theorem~\ref{genBDH}, we have the following theorem for square moduli.

\begin{theorem} \label{genBDHsquare}
For any $\varepsilon >0$ and fixed $A>0$, we have
\begin{equation} \label{genBDHsquareeq}
 \sum_{q \leq x^{1/2-\varepsilon}} q \sum_{\substack{a=1 \\ \gcd(a,q)=1}}^{q^2} \left| \psi(x,q^2,a) - \frac{x}{\varphi(q^2)} \right|^2 \ll \frac{x^2}{(\log x)^A}.
\end{equation}
\end{theorem}

Next, we have a Bombieri-Vinogradov type theorem for sparse sets of moduli.

\begin{theorem} \label{genbomvino}
Let $S$ be a well-distributed set.  Assume that condition \eqref{conditionforboth} holds and
\[ |S(Q)| \gg Q^{\frac{3}{4}} \; \mbox{as} \; Q \to \infty. \]
Then, if $Q \leq x^{\frac{1}{2}-\varepsilon}$, we have, for any fixed $A>0$,
\begin{equation} \label{genbomvinoeq}
 \sum_{q \in S(Q)} \max_{\substack{a \\ \gcd(a,q)=1}} \left| \psi( x; q,a) - \frac{x}{\varphi(q)} \right| \ll \frac{|S(Q)|}{Q} \frac{x}{(\log x)^A}.
\end{equation}

\end{theorem}

We shall also prove the following Bombieri-Vinogradov type theorem for square moduli.

\begin{theorem} \label{bomvinosquare}
For any $\varepsilon >0$ and fixed $A>0$, we have
\begin{equation} \label{bomvinosquareeq}
\sum_{q \leq x^{2/9-\varepsilon}} q \max_{\substack{a \\ \gcd(a,q)=1}}
\left| \psi(x;q^2,a) - \frac{x}{\varphi(q^2)} \right| \ll \frac{x}{(\log x)^A}.
\end{equation}
\end{theorem}

A better version of Theorem~\ref{bomvinosquare} can be obtained under the assumption of a conjecture, as will be noted in a later section. \newline

It should also be noted that we achieve a saving of an arbitrary power of logarithm beyond the trivial bound in our results, as in the classical results.  Moreover, the weight $q$ that is present in both \eqref{genBDHsquareeq} and \eqref{bomvinosquareeq} is the is the suitable for square moduli to obtain results analogous to the classical theorems. \newline

As an application of Theorem~\ref{bomvinosquare}, we have the following.  It is an old but hitherto unresolved problem if there exist infinitely many primes of the form $m^2+1$, where $m$ is an integer.  G. H. Hardy and J. E. Littlewood \cite{HaLi} gave a conjectural asymptotic estimate for the number of primes not exceeding $x$ of this form.  One may find several results on approximations to this problem in the literature.  For example, Ankeny \cite{Ank} and Kubilius \cite{Kub} showed independently that under the Riemann hypothesis for Hecke $L$-functions there exist infinitely many primes of the form $p=m^2+n^2$ with $n<c\log p$, where $c$ is some positive constant.  For some unconditional results in this direction see \cite{Bal} and \cite{HaLe}.  Moreover, the prime number theorem of Pyatecki\u\i-\v Sapiro \cite{PiaSa} can also be viewed as an approximation of the problem of representation of primes by quadratic polynomials. \newline

Here, we approximate the problem of representation of primes by $m^2+1$ in the following way.  For a natural number $n$ let $s(n)$ be the unique square-free number $a$ such that $\ratq(\sqrt{n})=\ratq(\sqrt{a})$. In other words, $s(n)=n/m^2$, where $m^2$ is the largest square dividing $n$. We note that $s(n)=1$ if and only if $n$ is a perfect square.  We consider primes of the form $n+1$, where $s(n)$ is small. More in particular, we shall prove the following.

\begin{theorem} \label{primes} 
Let $\varepsilon>0$. Then there exist infinitely many primes $p$ such that $s(p-1)\le p^{5/9+\varepsilon}$.
\end{theorem}

Equivalently, the above result may be formulated as follows. \newline

\noindent {\bf Theorem~\ref{primes}$^{\prime}$.} \begin{it} Let $\varepsilon>0$. Then there exist infinitely many primes of the form $p=am^2+1$ such that $a\le p^{5/9+\varepsilon}$. \end{it} \newline

We also note that the set of integers of the form $m^2+1$ is very sparse.  Friedlander and Iwaniec \cite{FrIw} proved the celebrated result that there exist infinitely many primes of the form $m^2+n^4$ (with an asymptotic formula). The set of integers of the form $m^2+n^4$ contains of the set of integers of the form $m^2+1$ but is still very sparse.  The number of such integers not exceeding $x$ is $O(x^{3/4})$.  It is noteworthy that the set of natural numbers $n$ with $s(n)\leq n^{5/9+\varepsilon}$ is also very sparse. More in particular, the number of $n\le x$ with $s(n)\le n^{5/9+\varepsilon}$ is $O(x^{7/9+\varepsilon/2})$ as the following calculation shows.
\[ |\{n\le x\ :\ s(n)\le n^{5/9+\varepsilon}\}| \le |\{(a,m)\in \natn^2\ :\ a\le x^{5/9+\varepsilon},\ am^2\le x\}| = \sum\limits_{a\le x^{5/9+\varepsilon}} \sum\limits_{m\le \sqrt{x/a}} 1 = O(x^{7/9+\varepsilon/2}). \]

\section{Preliminary Lemmas}

In this section, we quote lemmas that we shall need in the proofs of our theorems.  We begin with the following. \newline

\begin{lemma} [Siegel-Walfisz] \label{siewal}
Let $B$ be any positive constant.  Then there exist positive numbers $C_1(B)$ and $C_2(B)$ depending only on $N$, such that if $q$ satisfies
\[ q \leq (\log x)^B \]
then
\[ \psi(x; q,a) = \frac{x}{\varphi (q)} + O \left( x \exp \left( - C_1(B) \sqrt{\log x} \right) \right), \]
uniformly in $q$, and
\[ \psi(x, \chi) \ll x \exp \left( - C_2(B) \sqrt{\log x} \right), \; \mbox{where} \; \psi(x, \chi) = \sum_{n \leq x} \Lambda(n) \chi(n) \]
for any non principle character $\chi \pmod{q}$.
\end{lemma}

\begin{proof}
See section 22 of \cite{HD}.
\end{proof}

We shall also need the following identity.

\begin{lemma} [Vaughan Identity] \label{vaughanid}
Let $U \geq 1$, $V \geq 1$, $UV \leq x$.  Then for every arithmetic function $f(n)$, we have
\[ \sum_{n \leq x} f(n) \Lambda (n) \ll U (\log x) \max_{n \leq x} |f(n)| + (\log x)T_1(f) +T_2(f) + T_3(f) \]
with
\[ T_1(f) = \sum_{ l \leq \max(U,V)} \max_w \left| \sum_{w < k \leq x/l} f(kl) \right|, \; \mbox{and} \; T_i(f) = \left| \sum_{U < m \leq \max(x/V, UV)} \sum_{k \leq x/m} a_i (m) b_i (k) f(mk) \right| \]
for $ i=2 \mbox{ and } 3$, where $a_i(m)$, $b_i(k)$ are arithmetic functions which depend only on $U$ and $V$ and satisfy the inequalities $|b_i(k)| \leq \tau (k)$ and $|a_i(k)| \leq \log k$ for all $k \in \natn$.
\end{lemma}
\begin{proof}
This is Theorem 6.1.1 in \cite{Bru} and has its origin in \cite{RCV}.
\end{proof}

We shall also need the following lemma due to Polya and Vinogradov.

\begin{lemma} \label{polyvino}
Let $\chi(n)$ be a non-principal character modulo $q$ and $M \in \intz$ and $N \in \natn$.  Then we have
\[ \left| \sum_{M < n \leq M+N} \chi(n) \right| \leq 6 \sqrt{q} \log q. \]
\end{lemma}

\begin{proof}
This is Theorem 12.5 in \cite{HIEK}.
\end{proof}

We shall not need the classical large sieve inequality in its full generality as given by Davenport and Halberstam \cite{DH1}, but only the inequality touching upon Dirichlet characters.

\begin{lemma}[Classical Large Sieve] \label{classls}
Let $\{ a_n \}$ be a sequence of complex numbers.  Suppose that $M \in \intz$, $N, Q \in \natn$. Then we have
\begin{equation} \label{classlseq}
\sum_{q=1}^Q \frac{q}{\varphi(q)} \sideset{}{^{\star}}\sum_{\chi \bmod q} \left| \sum_{n=M+1}^{M+N} a_n \chi (n) \right|^2 \leq (Q^2 +N ) \sum_{n=M+1}^{M+N} |a_n|^2,
\end{equation}
where $\sideset{}{^{\star}}\sum$ henceforth denotes the sum over primitive characters to the specified modulus.
\end{lemma}
\begin{proof}
See \cite{MVa} for the proof of this lemma.
\end{proof}

We also need the following lemma.

\begin{lemma} \label{classest}
For $1\leq Q_0 < Q$ we have
\[ \sum_{Q_0 < q \leq Q} \frac{1}{\varphi(q)} \sideset{}{^{\star}}\sum_{\chi \bmod{q}} \max_{y \leq x} \left| \psi'(y, \chi) \right| \ll \left( \frac{x}{Q_0} + x^{\frac{5}{6}} \log Q + x^{\frac{1}{2}} Q \right) (\log Qx)^4. \]
\end{lemma}
\begin{proof}
See section 28 of \cite{HD}.
\end{proof}

We also need the following for the proof of Theorem~\ref{primes}.

\begin{lemma} \label{phi} We have
\begin{equation} \label{varphi}
\sum\limits_{y<q\le 2y} \frac{1}{\varphi(q^2)} = \frac{1}{2\zeta(2)y}+
O\left(\frac{\log y}{y^2}\right). 
\end{equation}
\end{lemma}
\vspace{0.5cm}

\begin{proof}
We write
\begin{eqnarray*}
& & \sum\limits_{y<q\le 2y} \frac{1}{\varphi(q^2)} = \sum\limits_{y<q\le 2y} \frac{1}{q\varphi(q)} = \sum\limits_{y<q\le 2y} \frac{1}{q^2}\sum\limits_{d\vert q} \frac{\mu(d)}{d} = \sum\limits_{d\le 2y} \frac{\mu(d)}{d} \sum\limits_{\substack{ y < q \leq 2y \\ q\equiv 0 \bmod{d}}} \frac{1}{q^2} \\
& = & \sum\limits_{d\le 2y} \frac{\mu(d)}{d} \sum\limits_{y/d<r\le 2y/d} 
\frac{1}{(dr)^2} = \sum\limits_{d\le 2y} \frac{\mu(d)}{d^3} \left(\frac{d}{2y}+ O\left(\frac{d^2}{y^2}\right)\right) 
= \frac{1}{2y} \sum\limits_{d=1}^{\infty} \frac{\mu(d)}{d^2} +
O\left(\frac{\log y}{y^2}\right).
\end{eqnarray*}
From this and
\[ \sum\limits_{d=1}^{\infty} \frac{\mu(d)}{d^2}=\frac{1}{\zeta(2)}, \]
we obtain \eqref{varphi}.
\end{proof}

\section{Large sieve inequality with sparse sets of moduli}

As stated before, large sieve is one of the key ideas in the proofs of our theorems.  We have the following lemma.

\begin{lemma} \label{lssparse}
Suppose that $S$ is a well-distributed set.  Let $\{ a_n \}$ be a sequence of complex numbers, $M \in \intz$, $Q, t \in \natn$ and $t \leq Q$.  Then we have
\begin{equation} \label{lssparseeq1}
\sum_{q\in S_t(Q/t)} \sum_{\substack{a=1 \\ \gcd(a,q)=1}}^q \left| \sum_{n=M+1}^{M+N} a_n e \left( \frac{a}{q} n \right) \right|^2 \ll \left( N + \frac{Q}{t} (QN)^{\varepsilon} \left( \sqrt{N}+|S_t(Q/t)| \right) \right) \sum_{n=M+1}^{M+N} |a_n|^2,
\end{equation}
where the implied constant depends at most on $\varepsilon$ and $S$.
\end{lemma}

\begin{proof}
We observe from \eqref{defwelldis} that $\left( S_{t_1} \right)_{t_2}=S_{t_1t_2}$.  Using this and \eqref{expected}, we deduce the lemma from Theorem 2 in \cite{Ba1} with $\mathcal{S}$ replaced by $S_t(Q/t)$ and $X=Q^{\varepsilon}$
after a linear change of variables in $n$.
\end{proof}

\begin{lemma} \label{lssparse2}
Under the same conditions as in Lemma~\ref{lssparse}, we have
\begin{equation} \label{lssparseeq2}
\sum_{q\in S_t(Q/t)} \frac{q}{\varphi(q)} \sideset{}{^{\star}}\sum_{\chi \bmod{q}} \left| \sum_{n=M+1}^{M+N} a_n \chi (n) \right|^2 \ll \left( N + \frac{Q}{t} (QN)^{\varepsilon} \left( \sqrt{N}+S_t(Q/t) \right) \right) \sum_{n=M+1}^{M+N} |a_n|^2,
\end{equation}
where the implied constant depends at most on $\varepsilon$ and $S$.
\end{lemma}

\begin{proof}
This lemma follows from Lemma~\ref{lssparse} in the usual way by considering Gauss sums.
\end{proof}

\begin{lemma} \label{lssparse3}
Let $\{ a_m \}$ and $\{ b_n \}$ be two sequences of complex numbers.  Then under the same conditions as Lemma~\ref{lssparse}, we have
\begin{equation} \label{lssparseeq3}
\sum_{q \in S_t(Q/t)} \frac{q}{\varphi(q)} \sideset{}{^{\star}}\sum_{\chi \bmod{q}} \max_X \left| \mathop{\sum_{m \leq M} \sum_{n \leq N}}_{mn \leq X} a_m b_n \chi (mn) \right| \ll \log (2MN) \left( \Delta(M, Q, t) \Delta(N, Q, t) Z_a Z_b \right)^{\frac{1}{2}} ,
\end{equation}
where the implied constant depends at most on $\varepsilon$ and $S$,
\[ \Delta ( Y, Q, t) =  Y + \frac{Q}{t} (QY)^{\varepsilon} \left( \sqrt{Y}+|S_t(Q/t)| \right), \; Z_a = \sum_{m \leq M} |a_m|^2 \; \mbox{and} \; Z_b = \sum_{n \leq N} |b_n|^2. \]
\end{lemma}

\begin{proof}
This lemma follows from Lemma~\ref{lssparse2} in the standard way, using the techniques in section 28 in \cite{HD}.
\end{proof}

\section{Proof of Theorem~\ref{genBDH}}

We are now ready to prove Theorem~\ref{genBDH}.  Following the arguments in section 29 of \cite{HD}, we have
\begin{equation}
\sum_{\substack{a=1 \\ \gcd(a,q)=1}}^q \left( \psi(x;q,a) - \frac{x}{\varphi(q)} \right)^2 \ll (\log qx)^2 + \frac{1}{\varphi(q)} \sum_{\chi} \left| \psi'\left( x, \chi_1 \right) \right|^2,
\end{equation}
where
\[ \psi' (y, \chi) = \left\{ \begin{array}{lc} \psi(y,\chi) & \mbox{if} \; \chi \neq \chi_0, \\ \psi(y,\chi_0) -y & \mbox{if} \; \chi = \chi_0, \end{array} \right. \; \mbox{and} \]
$\chi_1$ is the primitive character that induces $\chi$.  Therefore, we have
\begin{equation} \label{homer}
\sum_{q \in S(Q)} \sum_{\substack{a=1 \\ \gcd(a,q)=1}}^q \left( \psi(x;q,a) - \frac{x}{\varphi(q)} \right)^2 \ll \sum_{\substack{t,q \\ Q \leq tq \leq 2Q \\ tq \in S(Q)}} \frac{1}{\varphi(q)\varphi(t)} \sideset{}{^{\star}}\sum_{\chi \bmod{q}} \left| \psi' (x, \chi) \right|^2 + \sum_{q \in S(Q)} (\log qx)^2.
\end{equation}

The second term on the right-hand side of \eqref{homer} is
\begin{equation} \label{minorterm}
 \ll |S(Q)| (\log Qx)^2.
\end{equation}
It now remains to estimate the first term on the right-hand side of \eqref{homer}.  If $Q \leq (\log x)^B$ for a large fixed $B>0$, then we apply Siegel-Walfisz Theorem, Lemma~\ref{siewal}, for the estimate.  Otherwise, we divide the sum of interest into the following two sums.
\[ \Sigma_1 = \sum_{q \leq 2Q (\log x)^C/|S(Q)|} \frac{1}{\varphi(q)} \sum_{t \in S_q(Q/q)} \frac{1}{\varphi(t)} \sideset{}{^{\star}}\sum_{\chi \bmod{q}} \left| \psi' (x, \chi) \right|^2 \]
and
\[ \Sigma_2 = \sum_{t \leq |S(Q)|/(\log x)^C} \frac{1}{\varphi(t)} \sum_{q \in S_t(Q/t)} \frac{1}{\varphi(q)} \sideset{}{^{\star}}\sum_{\chi \bmod{q}} \left| \psi' (x, \chi) \right|^2, \]
where $2<C<B/2$.  Note that $\chi$ is never trivial in the inner-most sum of $\Sigma_2$ since $q>1$.  Using \eqref{conditionforboth}, we have
\begin{equation} \label{lisa}
 \sum_{t \in S_q(Q/q)} \frac{1}{\varphi(t)} \ll (\log \log Q) \left| S_q (Q/q) \right| \frac{q}{Q} \ll (\log \log Q) |S(Q)| \frac{q^{1-\varepsilon}}{Q}.
\end{equation}
Therefore, we have
\[ \Sigma_1 \ll \frac{(\log \log Q) |S(Q)|}{Q} \sum_{q \leq 2Q (\log x)^C/|S(Q)|} \frac{q^{1-\varepsilon}}{\varphi(q)} \sideset{}{^{\star}}\sum_{\chi \bmod{q}} \left| \psi' (x, \chi) \right|^2. \]
Using Siegel-Walfisz Theorem, Lemma~\ref{siewal}, for $q \leq (\log x)^D$ with an appropriate $D>0$ and the classical large sieve in Lemma~\ref{classls} for $q > (\log x)^D$, we obtain
\[ \Sigma_1 \ll \frac{|S(Q)|}{Q} \left( \frac{Q^2 (\log x)^{2C}}{|S(Q)|^2} + x \right) \frac{x}{(\log x)^A} \ll \frac{|S(Q)|}{Q} \frac{x^2}{(\log x)^A}, \]
where we have used \eqref{conditionfor1} and $Q\le x^{1-\varepsilon}$
for the second inequality above.  To estimate $\Sigma_2$, we use Lemma~\ref{lssparse2}, the fact that $|S_t(Q/t)| \leq |S(Q)|$ and the conditions of the 
theorem to get
\begin{eqnarray*}
\Sigma_2 & \ll & \sum_{t \leq |S(Q)|/(\log x)^C} \frac{t}{\varphi(t)Q} \left( x + \frac{Q}{t} (Qx)^{\varepsilon} \left( \sqrt{x} + |S_t(Q/t)| \right) \right) x \log x \\
 & \ll & \frac{|S(Q)|x^2}{Q(\log x)^A} + x^{3/2+\varepsilon}+x^{1+\varepsilon} |S(Q)| \ll \frac{|S(Q)|x^2}{Q(\log x)^A},
\end{eqnarray*}
where $A=C-2$.  Combining the estimates for $\Sigma_1$, $\Sigma_2$ and \eqref{minorterm}, we get the desired result.

\section{Proof of Theorem~\ref{genbomvino}}

We will now prove Theorem~\ref{genbomvino}.  Following the arguments in section 28 of \cite{HD}, we have
\begin{equation}
\max_{\substack{a \\ \gcd(a,q)=1}} \left| \psi(x; q,a) - \frac{x}{\varphi(q)} \right| \ll (\log qx)^2 + \frac{1}{\varphi(q)} \sum_{\chi \bmod{q}} \left| \psi'(x, \chi_1) \right|,
\end{equation}
where $\psi' (x, \chi)$ is as defined in the previous section and $\chi_1$ is the primitive character that induces $\chi$.  Therefore, we have 
\begin{equation} \label{marge}
 \sum_{q \in S(Q)} \max_{\substack{a \\ \gcd(a,q)=1}} \left| \psi( x; q,a) - \frac{x}{\varphi(q)} \right| \ll \sum_{\substack{t,q \\ Q \leq tq \leq 2Q \\ tq \in S(Q)}} \frac{1}{\varphi(q)\varphi(t)} \sideset{}{^{\star}}\sum_{\chi \bmod{q}} \left| \psi' (x, \chi) \right| + \sum_{q \in S(Q)} (\log qx)^2.
\end{equation}
The second term on the right-hand side of \eqref{marge} is
\begin{equation} \label{bart}
 \ll |S(Q)| (\log Qx)^2.
\end{equation}
It still remains to estimate the first term on the right-hand side of \eqref{marge}.  If $Q \leq (\log x)^B$ for a large fixed $B>0$, then we apply Siegel-Walfisz Theorem, Lemma~\ref{siewal}, for the estimate.  Otherwise, we divide the sum of interest into the following two sums.
\begin{equation} \label{udo}
 \Sigma_1 = \sum_{q \leq 2Q/Q_1} \frac{1}{\varphi(q)} \sum_{t \in S_q(Q/q)} \frac{1}{\varphi(t)} \sideset{}{^{\star}}\sum_{\chi \bmod{q}} \left| \psi' (x, \chi) \right| 
\end{equation}
and
\begin{equation} \label{otto}
 \Sigma_2 = \sum_{t \leq Q_1} \frac{1}{\varphi(t)} \sum_{q \in S_t(Q/t)} \frac{1}{\varphi(q)} \sideset{}{^{\star}}\sum_{\chi \bmod{q}} \left| \psi (x, \chi) \right|,
\end{equation}
where $Q_1<Q$ is to be chosen later.  Note again here that the inner most sum of $\Sigma_2$ never contains the trivial character.  By the virtue of \eqref{lisa}, we have
\[ \Sigma_1 \ll \frac{(\log \log Q) |S(Q)|}{Q} \sum_{q \leq 2Q/Q_1} \frac{q^{1-\varepsilon}}{\varphi(q)} \sideset{}{^{\star}}\sum_{\chi \bmod{q}} \left| \psi' (x, \chi) \right|. \]
Using Siegel-Walfisz Theorem, Lemma~\ref{siewal}, for $q \leq (\log x)^C$ with an appropriate $C>0$ and Lemma~\ref{classest} for $q > (\log x)^C$, we obtain
\begin{equation} \label{flanders}
 \Sigma_1 \ll \frac{|S(Q)|}{Q} \left( \frac{x}{(\log x)^A} + x^{\frac{5}{6}} \left( \frac{Q}{Q_1} \right)^{1-\varepsilon} + \sqrt{x}\left( \frac{Q}{Q_1} \right)^{2-\varepsilon} \right) (\log x Q)^5.
\end{equation}
To estimate $\Sigma_2$, we majorize $\Sigma_2$ in the following way.
\begin{equation} \label{grandpa}
\Sigma_2 \ll \frac{1}{Q} \sum_{t \leq Q_1} \frac{t}{\varphi(t)} \Sigma_2(t), \; \mbox{where}\; \Sigma_2(t) = \sum_{q \in S_t(Q/t)} \frac{q}{\varphi(q)} \sideset{}{^{\star}}\sum_{\chi \bmod{q}} \left| \psi(x, \chi) \right|.
\end{equation}
To estimate $\Sigma_2(t)$, we proceed similarly as in the proof of the classical Bombieri-Vinogradov theorem \eqref{classicbomvino}.  Using the Vaughan identity, Lemma~\ref{vaughanid}, we have
\begin{equation} \label{aftervaughanid}
\Sigma_2(t) \ll U (\log x) \frac{Q}{t} \left| S_t(Q/t) \right| + (\log x) \Sigma_{2,1} (t) + \Sigma_{2,2}(t) + \Sigma_{2,3}(t),
\end{equation}
where
\begin{equation} \label{typeisum}
\Sigma_{2,1} (t) = \sum_{q \in S_t(Q/t)} \frac{q}{\varphi(q)} \sideset{}{^{\star}}\sum_{\chi \bmod{q}} \sum_{l \leq \max(U,V)} \max_w \left| \sum_{w <k \leq x/l} \chi(kl) \right|
\end{equation}
and
\begin{equation} \label{typeiisum}
\Sigma_{2,i} (t) = \sum_{q \in S_t(Q/t)} \frac{q}{\varphi(q)} \sideset{}{^{\star}}\sum_{\chi \bmod{q}} \left| \sum_{U\leq m  \leq \max(UV,x/V)} \sum_{k \leq x/m} a_i(m) b_i(k) \chi(km) \right|,
\end{equation}
with $i=2,3$, $|a_i(m)| \leq \log x$ and $|b_i(k)| \leq \tau(k)$.  \eqref{typeisum} can be estimated using Polya-Vinogradov theorem, Lemma~\ref{polyvino}.  We have
\begin{equation} \label{estitypei}
\Sigma_{2,1} (t) \ll \max(U,V) \left| S_t(Q/t) \right| \left( \frac{Q}{t} \right)^{\frac{3}{2}} \log x.
\end{equation}
We split the sum over $m$ in \eqref{typeiisum} into dyadic intervals and estimate the following.
\begin{equation} \label{typeiisumsplit}
\Sigma_{2,i}(t, M) = \sum_{q \in S_t(Q/t)} \frac{q}{\varphi(q)} \sideset{}{^{\star}}\sum_{\chi \bmod{q}} \left| \sum_{M < m \leq 2M} \sum_{k \leq x/m} a_i(m) b_i(k) \chi(km) \right|.
\end{equation}
By the virtue of Lemma~\ref{lssparse3}, we get
\begin{equation} \label{maggie}
\Sigma_{2,i}(t, M) \ll (\log x)^4 \sqrt{x} \left( \frac{Q}{t} x^{\varepsilon} \left( |S_t(Q/t)| + \sqrt{M} \right) + M \right)^{\frac{1}{2}} \times \left( \frac{Q}{t} x^{\varepsilon} \left( |S_t(Q/t)| + \sqrt{\frac{x}{M}} \right) + \frac{x}{M} \right)^{\frac{1}{2}}.
\end{equation}

Expanding \eqref{maggie}, summing up the contributions from the dyadic intervals, setting $U=V$ and using \eqref{estitypei}, \eqref{aftervaughanid} and \eqref{grandpa}, we obtain
\begin{equation} \label{herb}
\begin{split}
\Sigma_2 \ll \; x^{\varepsilon} &Q^{\frac{1}{2}} U \Omega_{1,3/2} + x^{\frac{1}{2} + \varepsilon} \Omega_{1,1} + \left( \frac{x^{3/4+\varepsilon}}{U^{1/4}} + x^{\frac{1}{2}+\varepsilon} U^{\frac{1}{2}} \right) \Omega_{1/2,1} + x^{\frac{3}{4}+\varepsilon} \\
& + \left( \frac{x^{1/2+\varepsilon}U}{Q^{1/2}} + \frac{x^{1+\varepsilon}}{Q^{1/2}U^{1/2}} \right) \Omega_{1/2,1/2} + \left( \frac{x^{3/4+\varepsilon}U^{1/2}}{Q^{1/2}} + \frac{x^{1+\varepsilon}}{Q^{1/2}U^{1/4}} \right) Q_1^{1/2} + \frac{x(\log x)^5 Q_1}{Q},
\end{split}
\end{equation}
where
\[ \Omega_{a,b} = \sum_{t \leq Q_1} \frac{|S_t(Q/t)|^a}{t^b}. \]
Since $|S_t(Q/t)| \leq |S(Q)|$, we have
\[ \Omega_{1,3/2} \ll |S(Q)|, \; \Omega_{1,1} \ll |S(Q)| \log x, \; \Omega_{1/2,1} \ll |S(Q)|^{\frac{1}{2}} \log x. \]
Moreover, we have
\begin{equation*}
 \Omega_{1/2,1/2} \leq \left( \sum_{t \leq Q_1} \frac{1}{t} \right)^{\frac{1}{2}} \left( \sum_{t \leq Q_1} |S_t(Q/t)| \right)^{\frac{1}{2}} \ll (\log x)^{\frac{1}{2}} \left( \sum_{q \in S} \sum_{t|q} 1 \right)^{\frac{1}{2}} \ll x^{\varepsilon} |S(Q)|^{\frac{1}{2}}.
\end{equation*}
Consequently, the right-hand side of \eqref{herb} is
\begin{equation} \label{millhouse}
\begin{split}
\ll \; x^{\varepsilon} & Q^{\frac{1}{2}} U |S(Q)| + x^{\frac{1}{2} + \varepsilon} |S(Q)| + \left( \frac{x^{3/4+\varepsilon}}{U^{1/4}} + x^{\frac{1}{2}+\varepsilon} U^{\frac{1}{2}} \right) |S(Q)|^{\frac{1}{2}} + x^{\frac{3}{4}+\varepsilon} \\
& + \left( \frac{x^{1/2+\varepsilon}U}{Q^{1/2}} + \frac{x^{1+\varepsilon}}{Q^{1/2}U^{1/2}} \right) |S(Q)|^{\frac{1}{2}} + \left( \frac{x^{3/4+\varepsilon}U^{1/2}}{Q^{1/2}} + \frac{x^{1+\varepsilon}}{Q^{1/2}U^{1/4}} \right) Q_1^{1/2} + \frac{x(\log x)^5 Q_1}{Q}.
\end{split}
\end{equation}
Choosing $U=x^{1/4}$ and using \eqref{millhouse} and the conditions of the theorem, we get
\begin{equation} \label{mrburns}
\Sigma_2 \ll \frac{|S(Q)|x^{1-\varepsilon}}{Q} + \frac{x^{15/16+ \epsilon}}{Q^{1/2}} Q_1^{\frac{1}{2}} + \frac{x (\log x)^5 Q_1}{Q}.
\end{equation}
Now we choose
\[ Q_1 = \min \left\{ \frac{x^{1/8-3\varepsilon}|S(Q)|^2}{Q}, \; \frac{|S(Q)|}{(\log x)^{A+5}} \right\}. \]
Using \eqref{marge}, \eqref{bart}, \eqref{flanders}, \eqref{mrburns} and the conditions of the theorem, we get the desired result.

\section{Proof of Theorem~\ref{bomvinosquare}}

Let $S$ be the set of squares.  It suffices to estimate \eqref{udo} and \eqref{otto} with square moduli under the assumption $Q > (\log x)^B$ as it does in the previous section. \newline

Using the first inequality in \eqref{lisa} and the estimate
\[ |S_q(Q/q)| \ll \frac{|S(Q)|}{q^{1/2}} \ll \left( \frac{Q}{q} \right)^{\frac{1}{2}} , \]
which follows from section 6 of \cite{Ba1}, we obtain
\[ \Sigma_1 \ll \frac{\log \log Q}{Q^{1/2}} \sum_{q \leq 2Q/Q_1} \frac{q^{1/2}}{\varphi(q)} \sideset{}{^{\star}}\sum_{\chi \bmod{q}} \left| \psi'(x,\chi) \right|. \]
Using Siegel-Walfisz, Lemma~\ref{siewal}, for small moduli $q \leq (\log x)^B$ 
with an appropriate $B>0$ and Lemma~\ref{classest} for $q > (\log x)^B$, we obtain
\begin{equation} \label{jasper}
\Sigma_1 \ll \frac{(\log xQ)^5}{Q^{1/2}} \left( \frac{x}{(\log x)^A} + x^{\frac{5}{6}} \left( \frac{Q}{Q_1} \right)^{\frac{1}{2}} + x^{\frac{1}{2}} \left( \frac{Q}{Q_1} \right)^{\frac{3}{2}} \right).
\end{equation}
To estimate $\Sigma_2$, we use \eqref{millhouse}.  Choosing $U=x^{\frac{1}{3}}$ and using $|S(Q)| \ll Q^{\frac{1}{2}}$, we obtain
\begin{equation} \label{bea}
 \Sigma_2 \ll x^{\varepsilon} \left( Qx^{\frac{1}{3}} + Q^{\frac{1}{4}} x^{\frac{2}{3}} + x^{\frac{3}{4}} + x^{\frac{5}{6}} Q^{-\frac{1}{4}} + x^{\frac{1}{2}} Q^{\frac{1}{2}} + x^{\frac{11}{12}} Q^{-\frac{1}{2}} Q_1^{\frac{1}{2}} \right) + x (\log x)^5 Q_1Q^{-1}.
\end{equation}
Now we choose $Q_1 = Q^{\frac{3}{8}-\varepsilon}$.  Then if $(\log x)^C < Q \leq x^{\frac{4}{9}-\varepsilon}$ with a large $C>0$, we obtain
\[ \Sigma_1 + \Sigma_2 \ll \frac{x}{Q^{1/2} (\log x)^A} \]
from \eqref{jasper} and \eqref{bea}.  We have therefore proved, in a manner similar to the previous section,
\begin{equation} \label{willie}
\sum_{Q < q^2 \leq 2Q} \max_{\gcd(a,q)=1} \left| \psi(x;q^2,a) - \frac{x}{\varphi(q^2)} \right| \ll \frac{x}{Q^{1/2} (\log x)^A},
\end{equation}
if $(\log x)^C < Q \leq x^{\frac{4}{9}-\varepsilon}$ with $C>0$ chosen suitably.  \eqref{willie} also holds if $Q < (\log x)^C$ by Siegel-Walfisz Theorem, Lemma~\ref{siewal}.  \eqref{willie} implies now the desired estimate of the theorem.

\section{Proof of Theorem~\ref{primes}}

We shall derive Theorem~\ref{primes} from the following estimate.

\begin{lemma} \label{Lambda}
Let $\varepsilon$ and $A$ be arbitrary positive constants. If $(\log x)^{A+1}\le y\le x^{2/9-\varepsilon}$, then
\begin{equation} \label{asymp}
\sum\limits_{x<n\le 2x} \Lambda(n+1)\sum\limits_{\substack{y<q\le 2y\\ q^2|n}} 1 = \frac{x}{2\zeta(2)y} + O\left(\frac{x}{y(\log x)^A}\right),
\end{equation}
where the $O$-constant depends only on $\varepsilon$ and $A$. 
\end{lemma}   

\begin{proof}
To establish \eqref{asymp}, we write
\begin{eqnarray}
\nonumber & & \sum_{x<n\le 2x} \Lambda(n+1) \sum_{\substack{y<q\le 2y \\ q^2|n}} 1 = \sum_{y<q\le 2y} \;\sum_{\substack{x+1<k\le 2x+1 \\ k\equiv 1 \bmod{q^2}}} \Lambda(k) = \sum_{y<q\le 2y} (\psi(2x+1;q^2,1)-\psi(x+1;q^2,1)) \\\label{red} &=& x\sum\limits_{y<q\le 2y} \frac{1}{\varphi(q^2)} + O\left(\frac{1}{y}\sum_{q \leq 2y} q \left| \psi(2x+1;q^2,1) - \frac{2x+1}{\varphi(q^2)} \right| + \frac{1}{y}\sum_{q \leq 2y} q \left| \psi(x+1;q^2,1) - \frac{x+1}{\varphi(q^2)} \right|\right).
\end{eqnarray}
Now \eqref{asymp} follows from \eqref{red}, Theorem~\ref{bomvinosquare} and Lemma~\ref{phi}.
\end{proof}

Let $f$ be the characteristic function of the set of primes,
that is,
\[ f(n)=\left\{\begin{array}{llll} 1, & \mbox{ if } n \mbox{ is a prime,} \\ \\ 0, & \mbox{ otherwise.} \end{array} \right. \]
Furthermore, let $\tau(n)$ be the number of divisors of $n$. We have $s(n)\le n^{5/9+\varepsilon}$ if and only if there exists $q\in \natn$ such that $q^2|n$ and $q\ge n^{2/9-\varepsilon/2}$.
Therefore the number of primes $p\le x$ such that $s(p-1)\le p^{5/9+\varepsilon}$ is

\begin{equation} \label{trick}
\ge \sum\limits_{n\le x-1} f(n+1) \frac{1}{\tau(n)}\sum\limits_{\substack{q\ge n^{2/9-\varepsilon/2}\\ q^2|n}} 1.
\end{equation}
Using $\tau(n)\ll n^{\varepsilon/2}$ and the definition of the von Mangoldt function $\Lambda(n)$, the expression in (\ref{trick}) is 
\begin{equation} \label{start}
\ge \left(c_1 x^{-\varepsilon/2}\sum\limits_{n\le x-1} \Lambda(n+1) \sum\limits_{\substack{q\ge n^{2/9-\varepsilon/2}\\ q^2|n}} 1\right) - c_2x^{1/2+\varepsilon}
\end{equation}
for some positive constants $c_1,c_2$. Moreover, we have
\begin{equation} \label{that's it}
\sum\limits_{n\le x-1} \Lambda(n+1) \sum\limits_{\substack{q\ge n^{2/9-\varepsilon/2}\\ q^2|n}} 1 \ge \sum\limits_{(x-1)/2<n\le x-1} \Lambda(n+1)
\sum\limits_{\substack{y<q\le 2y\\ q^2|n}} 1,
\end{equation}
where $y:=x^{2/9-\varepsilon/2}$. From \eqref{trick}, \eqref{start}, \eqref{that's it} and Lemma~\ref{Lambda}, we deduce that the number of primes 
$p=n+1\le x$ with $s(n)\le n^{5/9+\varepsilon}$ is $\gg x^{7/9}$. This implies the result of  Theorem~\ref{primes}.

\section{Implication of a Conjecture}

We start with the following conjecture which can be viewed as a generalization of the conjecture given in \cite{Zha}.
\begin{conjecture} \label{conj}
Let $S$ be the set of squares and $S_t(Q/t)$ as defined before.  Let $M \in \intz$, $N \in \natn$ and $\{ a_n \}$ be a complex sequence.  Then we have
\begin{equation} \label{seymore}
\sum_{q \in S_t(Q/t)} \sum_{\substack{a=1 \\ \gcd(a,q)=1}}^q \left| \sum_{n =M+1}^{M+N} a_n e \left( \frac{a}{q} n \right) \right|^2 \ll Q^{\varepsilon} \left( \frac{Q}{t} \left| S_t \left( Q/t \right) \right| +N \right) \sum_{n=M+1}^{M+N} |a_n|^2,
\end{equation}
where the implied constant depends on $\varepsilon$.
\end{conjecture}

Under the truth of this conjecture, the range of summation for the outer-most sum in \eqref{bomvinosquareeq} can be extended to $q \leq x^{1/4-\varepsilon}$.  This better range would give the complete analogue of the classical Bombieri-Vinogradov theorem, up to an $\varepsilon$-factor, for square moduli. \newline

Furthermore, from \eqref{bomvinosquareeq} with the extended range for $q$ with $q \leq x^{1/4-\varepsilon}$, it would follow that there exist infinitely many primes $p$ such that $s(p-1)\le p^{1/2+\varepsilon}$.  We would get the same result under the assumption of the generalized Riemann hypothesis for Dirichlet $L$-functions.  We note that the set of $n$ such that $s(n)\le n^{1/2+\varepsilon}$ is ``almost'' as sparse as the set of numbers $m^2+n^4$ considered by Friedlander and Iwaniec \cite{FrIw}. Indeed, the number of $n\le x$ such that $s(n)\le n^{1/2+\varepsilon}$ is $O(x^{3/4+\varepsilon/2})$. \newline

It is conceivable that an Elliott-Halberstam type hypothesis holds for primes in arithmetic progressions to square moduli, {\it i.e.}, that \eqref{bomvinosquareeq} holds with the exponent $1/2-\varepsilon$ in place of $2/9-\varepsilon$. This would imply that there exist infinitely many primes $p$ such that $s(p-1)\le p^{\varepsilon}$.  A result of this kind comes very close to the conjecture that there exist infinitely many primes of the form $n^2+1$ since the number of $n\le x$ such that $s(n)\le n^{\varepsilon}$ is $O(x^{1/2+\varepsilon/2})$. \newline

{\bf Acknowledgement.}
This paper was written when the first and second-named authors held postdoctoral fellowships at the Department of Mathematics and Statistics at Queen's University and the Department of Mathematics at the University of Toronto, respectively.  The authors wish to thank these institutions for their financial support.

\bibliography{biblio}
\bibliographystyle{amsxport}

\vspace*{.7cm}

\noindent Department of Mathematics and Statistics, Queen's University \newline
99 University Ave, Kingston, ON K7L 3N6 Canada \newline
Email: {\tt sbaier@mast.queensu.ca} \newline

\noindent Department of Mathematics, University of Toronto \newline
40 Saint George Street, Toronto, ON M5S 2E4 Canada \newline
Email: {\tt lzhao@math.toronto.edu}
\end{document}